\documentclass[11pt, a4paper]{article}

\usepackage{latexsym}
\usepackage{graphicx, float}
\usepackage{amssymb,amsfonts,amsthm,amsmath,graphicx,url}

\long\def\delete#1{}

\newtheorem{theorem}{Theorem}

\newtheorem{proposition}[theorem]{Proposition}
\newtheorem{question}[theorem]{Question}

\input amssym.def
\input amssym.tex

\newcommand{\beq}{\begin{equation}}
\newcommand{\eeq}{\end{equation}}
\newcommand{\bea}{\begin{eqnarray}}
\newcommand{\eea}{\end{eqnarray}}
\newcommand{\bean}{\begin{eqnarray*}}
\newcommand{\eean}{\end{eqnarray*}}

\def\qed{\hfill$\Box$\vspace{10pt}}

\def\ZZZ{\Bbb Z}

\def\b0{{\bf 0}}

\def\De{\Delta}
\def\Ga{\Gamma}

\def\b{\beta}
\def\d{\delta}

\def\s{\sigma}
\def\t{\tau}
\def\ve{\varepsilon}

\def\PSL{{\rm PSL}}

\def\PG{{\rm PG}}

\def\PGammaL{{\rm P\Gamma L}}

\title{A note on the degree-diameter problem for arc-transitive graphs}
\author{Sanming Zhou\thanks{Supported by a Future Fellowship (FT110100629) of the Australian Research Council.}\\ 
Department of Mathematics and Statistics\\
The University of Melbourne\\
Parkville, VIC 3010, Australia\\
\small{\it smzhou@ms.unimelb.edu.au}}

\date{}

\begin{document}

\openup 0.5\jot

\maketitle

\begin{abstract}
We give two lower bounds on the largest order of an arc-transitive graph of diameter two and a given degree. 

\end{abstract}
  
All graphs in this note are finite and undirected with no loops or multiple edges. An {\em arc} of a graph is an ordered pair of adjacent vertices. A graph is {\em vertex-transitive} if its automorphism group is transitive on its vertex set, that is, any vertex can be permuted to any other vertex by an automorphism. A graph is {\em arc-transitive} if it is vertex-transitive and in addition its automorphism group is transitive on the set of arcs. 
  
Given integers $\De, D \ge 1$, let $N(\De, D)$ denote the maximum order (number of vertices) in a graph of maximum degree $\De$ and diameter at most $D$. The well known {\em degree-diameter problem} asks for finding $N(\De, D)$ together with the corresponding extremal graphs. A huge amount of work on this problem and its restrictions to several graph classes, including the classes of bipartite graphs, vertex-transitive graphs and Cayley graphs, has been done \cite{MS} since over fifty years ago \cite{HS}. In contrast, very little is known on this problem when restricted to the class of arc-transitive graphs \cite{wiki}. Denote by $N^{at}(\De, D)$ ($N^{vt}(\De, D)$, respectively) the maximum order of an arc-transitive graph (vertex-transitive graph, respectively) of degree $\De$ and diameter at most $D$. In this note we give two lower bounds on $N^{at}(\De, 2)$. 

In the case when $D=2$, the well known Moore bound \cite{MS} gives 
\beq
\label{eq:Moore}
N(\De, 2) \le \De^2 + 1
\eeq
for any $\De$. Here the equality holds only when $\De = 1, 2, 3, 7$ and possibly $57$, and for all other $\De$ we have $N(\De, 2) \le \De^2 - 1$ \cite{MS}. Regarding lower bounds, it is known \cite{Brown} that $N(\De, 2) \ge \De^2 - \De + 1$ for every $\De$ such that $\De - 1$ is a prime. In \cite{MMS} it was proved that, for $\De = (3q-1)/2$ with $q$ a prime power congruent to 1 modulo 4, we have
\beq
\label{eq:mms}
N^{vt}(\De, 2) \ge \frac{8}{9}\left(\De + \frac{1}{2}\right)^2.
\eeq
This bound came with the discovery \cite{MMS} of an infinite family of vertex-transitive graphs (now well known as the McKay-Miller-\v{S}ir\'{a}\v{n} graphs) with degree $\De = (3q-1)/2$ and order $8(\De + (1/2))^2/9$. Such extremal graphs are not arc-transitive except the Hoffman-Singleton graph. This is because, apart from this exceptional case, all McKay-Miller-\v{S}ir\'{a}\v{n} graphs admit a nontrivial partition \cite[Lemma 17]{Hafner} with cross-block edges which is invariant under the automorphism group, but on the other hand there are edges \cite[Definition 11]{Hafner} with end-vertices in the same block of the partition. Thus the lower bound in (\ref{eq:mms}) may not apply to $N^{at}(\De, 2)$. In \cite{ANS}, it was proved that for $\De = 2q-1$ with $q$ a prime power not congruent to 1 modulo 4, $N^{vt}(\De, 2) \ge (\De + 1)^2/2$.

Recently, a significant improvement over (\ref{eq:mms}) was made in \cite{SS}, where \v{S}iagiov\'{a} and \v{S}ir\'{a}\v{n} proved the existence of a Cayley graph of degree $\De = 2^{2m+\d}+(2+\d)2^{m+1} - 6$, diameter 2, and order greater than $\De^2 - 6\sqrt{2}\De^{3/2}$, where $m$ is any positive integer and $\d = 0$ or $1$. This implies, for $\De$ of this form, 
\beq
\label{eq:ss}
N^{vt}(\De, 2) \ge \De^2 - 6\sqrt{2}\De^{\frac{3}{2}}.
\eeq
Again this bound may not apply to $N^{at}(\De, 2)$ since the extremal graphs \cite{SS} are not arc-transitive. 

In view of (\ref{eq:Moore})-(\ref{eq:ss}), it is natural to ask \cite{Zhou-JGT} whether a lower bound on $N^{at}(\De, 2)$ of order $O(\De^2)$ exists for infinitely many $\De$ or even any $\De$. We record the following observations, the first of which answers this question in the affirmative.
\begin{proposition}
\label{main}
\begin{itemize}
\item[\rm (a)] For any even integer $\De \ge 2$,
\beq
\label{eq:hamming}
N^{at}(\De, 2) \ge \frac{1}{4}(\De + 2)^2.
\eeq
\item[\rm (b)] For any $0 < \ve < 1$, there exist infinitely many odd integers $\De$ of the form $q^{3}(q^{d-2} - 1)/(q-1)$, where $q$ is an odd prime power and $d \ge 3/\ve$ is an odd integer, such that 
\begin{equation}
\label{eq:p1}
N^{at}(\De, 2) > \De^{2 - \ve} + 2\De^{1 - \frac{\ve}{3}} + \De^{1 - \frac{2\ve}{3}} + 3. 
\end{equation}
\end{itemize}
\end{proposition} 

\thispagestyle{empty}

We remark that the arc-transitive graphs proving these results are not new. Nevertheless, it seems that these bounds have not been noted before; they are among the first general lower bounds for $N^{at}(\De, 2)$ as far as we know. It is our hope that these observations may inspire further studies on the degree-diameter problem for arc-transitive graphs.  

The graphs proving (\ref{eq:hamming}) are the well known {\em Hamming graph} $H(2, q)$ (where $q \ge 2$ is an integer), namely the graph with vertex set $\ZZZ_q^2$ such that two pairs are adjacent if and only if they differ in exactly one coordinate.  

The lower bound (\ref{eq:p1}) is obtained by using a result in \cite{Zhou-EJC}. 
Let $d \geq 3$ be an integer and $q$ a prime power. Two distinct lines $L, N$ of the projective space $\PG(d-1,q)$ are called {\em intersecting} if there exists a unique point incident with both $L$ and $N$ (that is, $L$ and $N$ lie on the same plane of $\PG(d-1,q)$).  
Define the {\em projective flag graph} $\Ga(d,q)$ to be the graph whose vertices are the (point, line)-flags of $\PG(d-1,q)$ such that two such flags $(\s, L), (\t, N)$ are adjacent if and only if $L, N$ are intersecting lines of $\PG(d-1,q)$. This graph was introduced in \cite[Definition 3.4]{Zhou-EJC} (and was denoted by $\Ga^{+}(P;d,q)$ there) in the classification of a family of arc-transitive graphs.  

\bigskip
\begin{proof}{\bf of Proposition \ref{main}:}~~
(a) Since $H(2, q)$ is arc-transitive (see e.g.~\cite{P}) with diameter 2, degree $\De=2(q-1)$ and order $q^2 = (\De/2 + 1)^2$, we obtain (\ref{eq:hamming}) immediately. Note that $\De \ge 2$ can be any even integer since $q \ge 2$ can be any integer. 

(b) Let $d \geq 3$ be an integer and $q$ a prime power. It is well known that any group $G$ with $\PSL(d,q) \leq G \leq \PGammaL(d,q)$ acts doubly transitively on the point-set of $\PG(d-1,q)$ (see e.g.~\cite{Beth-Jung-Lenz}). In \cite[Theorems 3.6-3.7]{Zhou-EJC} it was proved that the projective flag graph $\Ga(d,q)$ is $G$-arc transitive with diameter 2, girth 3 and degree $\De = q^{3}(q^{d-2} - 1)/(q-1)$. Since $\PG(d-1,q)$ has $(q^d - 1)/(q-1)$ points and each of them is in $(q^{d-1} - 1)/(q-1)$ lines (see e.g.~\cite[Proposition I.2.16]{Beth-Jung-Lenz}), the order of $\Ga(d,q)$ is equal to $N =  (q^d - 1)(q^{d-1} - 1)/(q-1)^2$. 

Note that $N = (q^{d-1}+ \cdots + q + 1)(q^{d-2}+ \cdots + q + 1) =  \{q^2 (q^{d-3}+ \cdots + q + 1)+(q+1)\}\{q(q^{d-3}+ \cdots + q + 1)+1\} = \left(\frac{\De}{q} + q+1\right) \left(\frac{\De}{q^2} + 1\right) = \frac{\De^2}{q^3} + \left(\frac{2}{q} + \frac{1}{q^2}\right) \De + (q+1)$. Since $\De = q^{3}(q^{d-2} - 1)/(q-1)$, we have $q \le \De^{1/d}$ and the equality holds if and only if $d=3$. Thus, for this $\De$, we have
$$
N^{at}(\De, 2) \ge \De^{2 - \frac{3}{d}} + 2\De^{1 - \frac{1}{d}} + \De^{1 - \frac{2}{d}} + 3 
$$
and the equality holds if and only if $d=3$ and $q=2$. 
Therefore, for any $0 < \ve < 1$, any prime power $q$, and any $d \ge 3/\ve$, we have 
$$
N^{at}(\De, 2) > \De^{2 - \ve} + 2\De^{1 - \frac{\ve}{3}} + \De^{1 - \frac{2\ve}{3}} + 3.
$$ 
In particular, if both $d$ and $q$ are odd, then $\De$ is odd and therefore (b) is proved. 
\qed
\end{proof}

\thispagestyle{empty}

Note that $\De$ is even when $d$ or $q$ is even. Although (\ref{eq:p1}) is still valid in this case, in general it is inferior to (\ref{eq:hamming}) in this case. 

In view of the discussion above we pose the following question.

\begin{question}
Are there infinitely many integers $\De \ge 2$ such that
$$
N^{at}(\De, 2) \ge \De^2 - f(\De)
$$
for some function $f$ with $f(x)/x^2 \rightarrow 0$ as $x \rightarrow \infty$?
\end{question}  

\medskip 
\subsection*{Acknowledgments}
Thanks go to G. Pineda-Villavicencio for helpful discussions on the degree-diameter problem, and to J. \v{S}ir\'{a}\v{n} for informing me of his recent work \cite{SS} with J. \v{S}iagiov\'{a}.  

\thispagestyle{empty}

\end{document}